\documentstyle[amstex,amssymb]{article} 
\newcommand{\ZZ}{{\mathbb Z}} 
\newcommand{\RR}{{\mathbb R}}  
\newcommand{\NN}{{\mathbb N}}  
\newcommand{\TT}{{\mathbb T}}    
\newcommand{\CC}{{\mathbb C}}    
\newtheorem{theorem}{Theorem}      
\newtheorem{lemma}{Lemma}[section]      
\newtheorem{prop}[lemma]{Proposition}      
      
\newtheorem{coro}[lemma]{Corollary}      
\begin{document}      
\title{Half-line eigenfunction estimates and  singular      
  continuous spectrum of zero Lebesgue measure}      
\author{David Damanik$\,^{1,2}$ and Daniel Lenz$\,^{2}$}      
\maketitle    
\begin{center}    
{\it Dedicated to Joachim Weidmann on the occasion of his 60th    
  birthday}    
\end{center}    
\vspace{0.5cm}    
${{}^1}$ Department of Mathematics,      
California Institute of Technology,      
Pasadena, CA 91125,      
U.S.A.\\[0.2cm]      
$^2$ Fachbereich Mathematik,      
Johann Wolfgang Goethe-Universit\"at,      
60054 Frankfurt, Germany\\[0.3cm]      
1991 AMS Subject Classification: 81Q10, 47B80\\      
Key words: Schr\"odinger operators, singular continuous spectrum, Cantor spectrum      
    
\begin{abstract}      
We consider discrete one-dimensional Schr\"odinger operators with      
strictly ergodic, aperiodic potentials taking finitely many values. The      
well-known tendency of these operators to have purely singular      
continuous spectrum of zero Lebesgue measure is further elucidated. We provide a unified approach to both the study of the spectral type as well as the measure of the spectrum as a set.  We apply this approach to Schr\"odinger operators with Sturmian potentials. Finally, in the appendix, we discuss the two different strictly ergodic dynamical systems associated to a circle map.     
\end{abstract}      
      
\section{Introduction}      
Consider a separable metric space $\Omega$, a Borel measure $\mu$ on      
$\Omega$ and a homeomorphism $T: \Omega \rightarrow \Omega$ such that      
$(\Omega, \mu, T)$ is ergodic. If $g:\Omega \rightarrow      
\RR$ is measurable, we define $V_\omega:\ZZ \rightarrow \RR$ by      
$V_\omega(n) = g(T^n \omega)$, $\omega \in \Omega$. This gives rise      
to an ergodic family of discrete one-dimensional Schr\"odinger operators      
$(H_\omega)_{\omega \in \Omega}$ in $l^2(\ZZ)$, namely,      
       
\begin{equation}\label{opfam}      
(H_\omega \phi)(n) = \phi (n+1) + \phi (n-1) + V_\omega (n) \phi (n).      
\end{equation}      
It is well known that the spectral properties of the operators      
$H_\omega$ are deterministic in the sense that there exist closed sets      
$\Sigma, \Sigma_{\rm pp}, \Sigma_{\rm sc}, \Sigma_{\rm ac} \subseteq \RR$ and a set      
$\Omega_0 \subseteq \Omega$ obeying $\mu(\Omega_0)=1$ such that the equations      
$\sigma (H_\omega) = \Sigma, \sigma_{\rm pp} (H_\omega) =      
\Sigma_{\rm pp},\sigma_{\rm sc} (H_\omega) = \Sigma_{\rm sc},\sigma_{\rm ac} (H_\omega) =      
\Sigma_{\rm ac}$ hold for every $\omega \in \Omega_0$. We will restrict our      
attention to the case where the dynamical system $(\Omega, \mu, T)$ is      
strictly ergodic and the resulting potentials $V_\omega$ are aperiodic and take only finitely many values. In this case, it follows from Kotani \cite{k2} that the set $\Sigma_{\rm ac}$ is actually empty,      
      
\begin{equation}      
  \label{acempty}      
  \sigma_{\rm ac} (H_\omega) = \emptyset \;\;\; \text{for a.e.~$\omega \in \Omega$}.      
\end{equation}      
Further investigations have shown that, in many cases, one also has empty      
point spectrum. On the other hand, no example with non-empty point      
spectrum is known. Moreover, in many cases the spectrum $\Sigma$ has      
zero Lebesgue measure, and again, no example is known where the spectrum      
is a set of non-vanishing Lebesgue measure. Thus, there is an apparent      
tendency for the operators $H_\omega$ to have purely singular continuous      
zero-measure spectrum. By general absence of absolutely continuous      
spectrum, the two basic properties ${\cal P}_1, {\cal P}_2$ one seeks to      
establish are therefore the following,      
$$\begin{array}[t]{ll}      
{\cal P}_1\; : & \Sigma_{\rm pp} = \emptyset,\\      
{\cal P}_2\; : & |\Sigma | = 0.      
\end{array}$$      
Of course, ${\cal P}_1$ can be weakened (resp., strengthened) to ${\cal      
  P}_1': \sigma_{\rm pp}(H_{\omega_0})=\emptyset$ for {\it some}       
$\omega_0 \in \Omega$ (resp., ${\cal P}_1'': \sigma_{\rm pp}(H_{\omega})=\emptyset$  for {\it every} $\omega \in \Omega$). We shall refer to ${\cal P}_1,      
{\cal P}_1', {\cal P}_1''$ as {\it almost sure, generic} and {\it      
  uniform} absence of eigenvalues, respectively. There are two important      
subclasses of families $(H_\omega)_{\omega \in \Omega}$ which we shall      
now briefly recall.\\[5mm]      
{\bf Potentials generated by primitive substitutions}: Let $A      
=\{a_1,\ldots,a_s\}$ be a finite set, called the {\it alphabet}. The      
$a_i$ are called {\it symbols} or {\it letters}, the elements of $A^*      
= \bigcup_{k \ge 1} A^k$ are called {\it words}. We denote by $|v|$      
the {\it length} of a word $v \in A^*$ (i.e., $|v|=l$ iff $v \in      
A^l$). Two words $w_1,w_2$ are called {\it conjugate}, denoted by $w_1 \sim w_2$, if they have the same length and $w_1$ can be transformed to $w_2$ by a cyclic permutation of its symbols. Let $A^\NN,A^\ZZ$ denote the set of one-sided and two-sided infinite sequences over $A$, respectively. A {\it substitution} $S$ is a map $S :       
A \rightarrow A^*$. $S$ can be extended homomorphically to $A^*$ (resp.,      
$A^{\NN}$) by $S(b_1 \ldots b_n) = S(b_1) \ldots S(b_n)$ (resp.,      
$S(b_1 b_2 b_3\ldots) = S(b_1) S(b_2) S(b_3)\ldots$). $S$ is called      
{\it primitive} if there exists $k \in \NN$ such that for every      
$a \in A$, $S^k(a)$ contains every symbol from $A$. Prominent examples      
of primitive substitutions are given by      
\begin{center}      
\begin{tabular}{ll}      
$a \mapsto ab$, $b \mapsto a$ & Fibonacci,\\      
$a \mapsto ab$, $b \mapsto aa$ & period doubling,\\      
$a \mapsto ab$, $b \mapsto aaa$ & binary non-Pisot,\\      
$a \mapsto ab$, $b \mapsto ba$ & Thue-Morse,\\      
$a \mapsto ab$, $b \mapsto ac$, $c \mapsto db$, $d \mapsto dc$ &      
Rudin-Shapiro.       
\end{tabular}      
\end{center}      
A fixed point $u \in A^{\NN}$ of $S$ is called {\it substitution      
  sequence}. The existence of such a fixed point is ensured by the      
following conditions,        
\begin{itemize}      
\item there exists a letter $a \in A$ such that the first letter of      
  $S(a)$ is $a$,      
\item $\lim_{n \rightarrow \infty} |S^n(a)| = \infty$,      
\end{itemize}      
which are easily seen to hold for a suitable power of $S$ if $S$ is      
primitive. Without loss of generality, we assume this power to be equal      
to one. In this case, $u = \lim_{n \rightarrow \infty} S^n(a)$      
exists and is a substitution sequence. Put discrete topology on      
$A$. Extend $u$ to the left arbitrarily giving $\hat{u} \in A^{\ZZ}$ and      
define $\Omega$ to be the set of accumulation points of the translates      
of $\hat{u}$ with respect to pointwise convergence,       
      
\begin{equation}\label{osubst}      
\Omega = \{ \omega \in A^{\ZZ} \; : \; \omega = \lim T^{n_i}      
\hat{u}, \; n_i \rightarrow \infty\},      
\end{equation}      
where the shift $T$ on $A^{\ZZ}$ is defined by $(T \tau)_n =\tau_{n+1}$,      
$\tau \in A^{\ZZ}$. By definition, $\Omega$ is a closed (and therefore      
compact) subset of $A^{\ZZ}$ that is invariant under $T$. The dynamical      
system $(\Omega , T)$ is called the {\it substitution dynamical system      
  associated to} $S$. If $S$ is primitive, this system is strictly      
ergodic \cite{q}. Recall that the Borel $\sigma$-algebra is generated by the cylinder sets  $[b_0 \ldots b_{l-1}]_{[m,m+l-1]} = \{ \omega \in \Omega : w_{m+i} = b_i, \; 0 \le i \le l-1\}$, $m \in \ZZ$, $l \ge 1$, $b_i \in A$, $0 \le i \le l-1$. Let us note that the unique ergodic Borel measure $\mu$ on $\Omega$ obeys       
      
\begin{equation}\label{muprop}      
\mu([b_0 \ldots b_{l-1}]_{[m,m+l-1]}) = d(b_0 \ldots b_{l-1}),      
\end{equation}      
where $d(b_0 \ldots b_{l-1})$ is the frequency of the word $b_0 \ldots b_{l-1}$ in $u$, that is,      
      
\begin{equation}      
d(b_0 \ldots b_{l-1}) = \lim_{n \rightarrow \infty} \frac{1}{n} \, \# \{j \le n \, : \, u_j \ldots u_{j + l - 1} = b_0 \ldots b_{l-1}\},        
\end{equation}      
which is always strictly positive. Pick an arbitrary function $f : A \rightarrow \RR$ and      
define $g:\Omega \rightarrow \RR$ by $g(\omega) = f(\omega_0)$. This      
yields potentials $V_\omega (n) = f(\omega_n)$.\\[5mm]       
{\bf Circle map potentials}: Let $\alpha \in (0,1)$ be irrational and      
consider on $\Omega = \TT = \RR / \ZZ \cong [0,1)$ the irrational      
rotation $T:\Omega \rightarrow \Omega,\; \omega \mapsto \omega +      
\alpha \mod 1$. The transformation $T$ is strictly ergodic with the Lebesgue      
measure on $\TT$ as unique ergodic measure $\mu$. Let $g$ be      
given by $g(\omega) = \lambda \chi_{[1-\beta,1)}(\omega)$, where $\beta \in (0,1)$. This yields potentials $V_\omega (n) = \lambda      
\chi_{[1-\beta,1)}(\alpha n+\omega \mod 1)$. There is another possibility to associate a family of potentials to the parameters $\alpha, \beta$ which is similar to the construction in the primitive substitution case. Namely, consider the function $v_{\alpha,\beta}(n)=\chi_{[1-\beta,1)}(n\alpha  \mod  1)$ and define the hull $\Omega_{\alpha,\beta}$ by      
      
\begin{equation}\label{ocircle}      
\Omega_{\alpha,\beta} = \{ \omega \in \{0,1\}^{\ZZ} \; : \; \omega = \lim T^{n_i} v_{\alpha,\beta}, \; n_i \rightarrow \infty\}.      
\end{equation}      
Now, the dynamical system $(\Omega_{\alpha,\beta},T)$ is strictly ergodic with unique ergodic measure $\mu$ given on cylinder sets by the frequency of the respective word, similar to the above. The function $g$ generating the potentials is in this case given by $g(\omega) = f(\omega_0)$, where $f(0)=0$, $f(1)=\lambda$. Both dynamical systems have been studied in the past, their mutual relation will be discussed in the last section of this paper. It will turn out that the two systems, while not being the same, are essentially equivalent, and we shall use the accumulation point representation in our study. The resulting potentials are called {\it Sturmian} if $\alpha = \beta$ and we shall in this case write $\Omega_\alpha$ instead of $\Omega_{\alpha,\alpha}$ and $v_\alpha$ instead of $v_{\alpha,\alpha}$. The case $\alpha = \frac{\sqrt{5}-1}{2}$ recovers the hull generated by the Fibonacci substitution. Schr\"odinger operators with Sturmian potentials provide a standard one-dimensional quasicrystal model \cite{lp1} and are therefore, from a physical point of view, the most interesting among the families generated by circle maps. See \cite{sbgc} for the report of the discovery of quasicrystals and \cite{s} for general background.\\[0.5cm]    
Let us remark that for families generated by primitive substitutions or circle maps, the potentials $V_\omega$ are minimal (i.e., for every pair $\omega_1,\omega_2 \in \Omega$, $V_{\omega_1}$ is a pointwise limit of translates of $V_{\omega_2}$). This implies the following (the claim for the spectrum is part of the folklore, and the part for the absolutely continuous spectrum follows from Last-Simon \cite{ls}),  
  
\begin{equation}\label{unifscm}  
\sigma (H_\omega) = \Sigma, \; \sigma_{\rm ac} (H_\omega) = \Sigma_{{\rm ac}} = \emptyset \;\;\; \text{for all $\omega \in \Omega$}.      
\end{equation}      
\\[5mm]  
Let us now summarize what is known in regard to ${\cal P}_1$ and ${\cal P}_2$ for the two classes introduced above. In the primitive substitution case, ${\cal P}_1$ was shown to be satisfied provided that there is a word $v \in A^*$ such that $vvvv$ occurs in the substitution sequence $u$ \cite{d3}. The argument applies to several prominent cases including period doubling (although a fourth power does not occur in this case) and binary non-Pisot. Generic absence of eigenvalues was shown for several examples \cite{bbg1,dp2} and subclasses \cite{bg1,bg2,hks}. Regarding ${\cal P}_2$, Bovier and Ghez have given a sufficient condition that is obeyed by a large subclass of primitive substitutions containing almost all the prominent examples, the Rudin-Shapiro substitution being a notable exception \cite{bg1}.\\[5mm]      
For circle map potentials, it follows from Delyon-Petritis \cite{dp1}      
and Kaminaga \cite{k1} that ${\cal P}_1$ holds for almost every $\alpha$      
and $\beta$ arbitrary and for all Sturmian families. Results of generic type were obtained in \cite{d1,hks,s3}. In particular, eigenvalues are always generically absent \cite{hks}. Recently, methods were developed that allow one to establish uniform absence of eigenvalues for all Sturmian hulls \cite{dl1,dl4}. On the other hand, Bellissard et al. have shown that ${\cal P}_2$ holds in the Sturmian case \cite{bist}, but no result is known for $\beta \not= \alpha$.\\[5mm]      
To summarize, it is by now quite well understood how to establish     
${\cal P}_1$, ${\cal P}_2$ in the primitive substitution case and      
${\cal P}_1$ in the circle map case, whereas it is not clear how to      
prove ${\cal P}_2$ for non-Sturmian circle map potentials. Now, our main objectives here are      
      
\begin{itemize}      
\item to propose a unified approach to both ${\cal P}_1$ and ${\cal      
P}_2$. In fact, we shall demonstrate that a certain strengthening of      
${\cal P}_1$ implies ${\cal P}_2$. This is particularly nice since there      
is a reasonable understanding of how to study ${\cal P}_1$ in the circle      
map case. To be precise, the shift of viewpoint will be the following:      
In order to exclude eigenvalues, one seeks to establish uniform {\it      
  lower} bounds for the solutions of the associated difference equation,      
whereas one usually looks for uniform {\it upper} bounds in order to      
prove ${\cal P}_2$. We will employ a theorem due to Osceledec which      
relates the two concepts and, in particular, reduces, in our context, the      
study of upper bounds to a study of lower bounds.      
\item to apply this unified approach to Schr\"odinger operators with    
  Sturmian potentials.      
\item to give a discussion of the two different strictly ergodic dynamical systems associated to a circle map.       
\end{itemize}      
The organization of this article is as follows. In Section 2 we      
define the stability set and show that if this set has positive measure,      
then both ${\cal P}_1$ and ${\cal P}_2$ hold. In Section 3 it is shown that this condition is satisfied for all Sturmian potentials. Section 4 provides extensions, remarks, and open problems, while Section 5 deals with a study of dynamical systems associated to circle maps.

\section{The stability set, singular continuous spectral measures, and      
  zero-measure spectrum}      
Consider a family $(H_\omega)_{\omega \in \Omega}$ such that the underlying dynamical system $(\Omega,T)$ is strictly egodic and the potentials $V_\omega$ are aperiodic and take finitely many values. Given the absence of absolutely continuous spectrum (\ref{acempty}), in order to prove purely singular continuous spectrum for $H_{\omega_0}$ it is sufficient to establish absence of $l^2$-solutions to the difference equation     
    
\begin{equation}      
  \label{eve}      
  u(n+1) + u(n-1) + V_{\omega_0} (n) u(n) = E u(n)      
\end{equation}      
for every $E \in \sigma(H_{\omega_0})$. Define the following stability sets:    
\begin{center}      
$\Omega^s_+$ $=$ $\{ \omega \in \Omega$ : for every $E \in \Sigma$,      
no solution of $(H_\omega - E)\phi = 0$ is square-summable at      
$+\infty\}$,      
\end{center}      
    
\begin{center}      
$\Omega^s_-$ $=$ $\{ \omega \in \Omega$ : for every $E \in \Sigma$,      
no solution of $(H_\omega - E)\phi = 0$ is square-summable at      
$-\infty\}$,      
\end{center}      
    
\begin{center}      
$\Omega^s$ $=$ $\{ \omega \in \Omega$ : $\forall E \in \Sigma \;\;      
\exists \,\varepsilon \in \{+,-\}$ s.t.       
no solution of $(H_\omega - E)\phi = 0$ is square-summable at      
$\varepsilon \infty\}$.      
\end{center}      
Obviously, we have the following inclusions,    
    
\begin{equation}      
  \label{omegas}      
  \Omega^s_+,\Omega^s_- \subseteq \Omega^s.      
\end{equation}      
The sets $\Omega^s_-, \Omega^s_+$ and $\Omega^s$ are called stability sets    
for the following reason: Consider an $\omega \in \Omega^s$ and let    
$H_\omega + W$ be a local perturbation of $H_\omega$, that is, $W$ is a   
finitely supported function $W: \ZZ \rightarrow \RR$. Then, an $E\in   
\Sigma$ can not be an eigenvalue of the $H_\omega +W$ . In this sense   
the absence of eigenvalues (on $\Sigma$) is stable under local   
perturbations for  $H_\omega$ with $\omega \in \Omega^s$. But, of   
course, a  finitely supported perturbation will in general introduce   
isolated eigenvalues outside of $\Sigma$ (compare the discussion at the end of this section). Let us summarize some basic properties which are immediate from the definitions.    
      
\begin{prop}      
For every $\omega \in \Omega^s$, no local perturbation of $H_\omega$ has eigenvalues in $\Sigma$.     
\end{prop}      
This applies in particular to $H_\omega$, $\omega \in \Omega^s$, itself. Naturally, we will be interested in the $\mu$-measure of the set      
$\Omega^s$. We shall show that the positivity of $\mu(\Omega^s)$ yields almost sure absence of eigenvalues, which is of course obvious, as well as zero-measure spectrum, and thus presents a unified treatment of the two basic problems of proving ${\cal P}_1$ and ${\cal P}_2$.

\begin{theorem}\label{basic}      
Suppose $\mu (\Omega^s ) > 0$. Then the family  $(H_\omega)_{\omega \in \Omega}$ satisfies ${\cal P}_1$ and ${\cal P}_2$. In particular, for  $\mu$-a.e. $\omega \in \Omega$, the operator $H_\omega$ has purely singular continuous spectrum of zero Lebesgue measure.     
\end{theorem}      
{\it Remark.} The condition $\mu (\Omega^s ) > 0$ means that there is a      
subset $\Omega_0 \subseteq \Omega^s$ which is measurable and has      
positive measure. Likewise, we shall say that $\mu (\Omega^s ) = 1$ if      
there is a measurable $\Omega_0 \subseteq \Omega^s$ having full      
measure.\\[5mm]      
Let us recall the definition of the Lyapunov exponent which is a central      
quantity in the study of one-dimensional ergodic Schr\"odinger      
operators. Given a family (\ref{opfam}) we define, for $E \in \CC$, $\omega      
\in \Omega$, $k,n \in \ZZ$, $k \le n$, transfer matrices $M(E,\omega,k,n)$ by       
$$M(E,\omega,k,n) =       
\left(       
  \begin{array}{cc}      
E-V_\omega(n) & -1\\1 & 0      
  \end{array}      
\right)      
\times      
\left(       
  \begin{array}{cc}      
E-V_\omega(n-1) & -1\\1 & 0      
  \end{array}      
\right)      
\times \cdots \times      
\left(       
  \begin{array}{cc}      
E-V_\omega(k) & -1\\1 & 0      
  \end{array}      
\right).$$      
The right/left-Lyapunov exponents $\gamma^\pm_\omega(E)$ are given by the      
following limits, provided that they exist,      
      
\begin{eqnarray*}      
\gamma^+_\omega(E) & = & \lim_{n \rightarrow \infty} \frac{1}{n} \ln \|M(E,\omega,1,n)\|,\\      
\gamma^-_\omega(E) & = & \lim_{n \rightarrow \infty} \frac{1}{n} \ln \|M(E,\omega,-n,-1)\|.      
\end{eqnarray*}      
It is well known that for every $E \in \CC$, there is a full measure      
set $\Omega_E \subseteq \Omega$ and a number $\gamma(E) \in \RR$ such      
that for every $\omega \in \Omega_E$, $\gamma^\pm_\omega(E)$ exist and      
the following equality holds,       
      
\begin{equation}\label{lyeq}      
\gamma^+_\omega(E) = \gamma^-_\omega(E) = \gamma(E)      
\end{equation}      
(see, e.g., \cite{cl}). The number $\gamma(E)$ is called {\it Lyapunov exponent}.      
\\[0.5cm]      
It will be seen from the proof of Theorem \ref{basic} that in order to      
establish ${\cal P}_2$, it suffices to find, for every $E\in \Sigma$, an      
$\omega$ such that $\omega$ belongs to the stability set and the      
right/left-Lyapunov exponents $\gamma^\pm_\omega(E)$ exist. Thus, the      
following theorem can be seen as a corollary to the proof of Theorem      
\ref{basic}.       
      
\begin{theorem}\label{const}      
Let $(H_\omega)_{\omega \in \Omega}$ be such that for every $E \in      
\Sigma$, $\Omega_E = \Omega$. Then, $\Omega^s \not= \emptyset$ implies ${\cal P}_2$.       
\end{theorem}      
{\it Remark.} The assumption of Theorem \ref{const} does not hold for      
general ergodic families. However, the theorem is of interest in our      
context since the assumption of uniform existence of the Lyapunov      
exponent for energies from the spectrum is obeyed by all families $(H_\omega)_{\omega \in \Omega}$ arising from a primitive substitution \cite{h} or from a Sturmian sequence \cite{dl2} and it is often relatively simple to find one element in $\Omega^s$, compare in particular \cite{bg1}. It is not known whether uniform existence of the Lyapunov exponent also holds for general circle map potentials.    
\\[0.5cm]      
{\it Proof of Theorem \ref{basic}.} Let $\Omega_0 \subseteq \Omega^s$ be      
measurable with $\mu(\Omega_0)>0$. Define $\Omega_1$ by       
$$\Omega_1 = \bigcup_{n \in \ZZ} T^n \Omega_0.$$      
Obviously, $\Omega_1$ is measurable, $T$-invariant, and has positive      
measure. Thus, $\mu(\Omega_1)=1$. On the other hand, it is also obvious      
that for every $n \in \ZZ$, the set $T^n \Omega_0$ is contained in      
$\Omega^s$. Thus, $\Omega_1 \subseteq \Omega^s$ and ${\cal P}_1$ holds. In order to prove ${\cal P}_2$, we recall the fundamental Kotani result $|\{E\;:\; \gamma(E) = 0\}| = 0$ \cite{k2}, $|\cdot|$ denoting Lebesgue measure. It thus suffices to prove       
    
\begin{equation}\label{sina}      
\Sigma \subseteq \{E\,:\, \gamma(E) = 0\}.      
\end{equation}      
Assume there exists an energy $E \in \Sigma$ such that $\gamma(E) >      
0$. Then, for a full measure set $\Omega_E \subseteq \Omega$,      
$\gamma^\pm_\omega(E)$ exist and obey (\ref{lyeq}). By $\mu(\Omega^s) >      
0$, the intersection $\Omega_E \cap \Omega^s$ is non-empty. Pick      
$\omega_0 \in \Omega_E \cap \Omega^s$ arbitrarily. Consider the operator      
$H_{\omega_0}$ and the corresponding eigenvalue equation      
      
\begin{equation}      
  \label{evesuit}      
  u(n+1) + u(n-1) + V_{\omega_0} (n) u(n) = E u(n).      
\end{equation}      
By $\omega_0 \in \Omega^s$, there exists a suitable choice of a half-line      
such that (\ref{evesuit}) has no solution which is square-summable on      
this half-line. We may assume w.l.o.g. that there are no solutions which      
are square-summable at $+\infty$. In this case we obtain a contradiction      
by $\omega \in \Omega_E$, (\ref{lyeq}), and Osceledec's theorem which      
states that $\gamma^+_{\omega_0}(E) > 0$ implies the existence of a      
solution of (\ref{evesuit}) which decays exponentially at $+\infty$      
\cite{cl} (see also \cite{ls} for a generalization).\hfill$\Box$\\[0.5cm]      
{\it Proof of Theorem \ref{const}.} In case $\Omega_E = \Omega$ for      
every $E \in \Sigma$, the condition $\Omega^s \not= \emptyset$ is      
sufficient to guarantee the existence of $\omega_0 \in \Omega_E \cap      
\Omega^s$ for each $E \in \Sigma$. For a proof of ${\cal P}_2$ we may      
then proceed as in the general case.\hfill$\Box$\\[0.5cm]      
We have seen that a study of solutions on a half-line enables us to      
establish the existence of purely singular continuous spectrum, and to provide a unified approach to both the spectral type and the Lebesgue measure of the spectrum as a set. We are now going to provide the fundamental strategy for studying the stability set      
$\Omega^s$. Basically, there are three approaches to deduce      
eigenfunction estimates from local properties of the potential: the      
two-block method, the three-block method (both being variants of an      
argument originally due to Gordon \cite{g}), and the palindrome method      
(due to Jitomirskaya-Simon \cite{js} and Hof et al.      
\cite{hks}). While we postpone a discussion and comparison of these      
methods to a subsequent section, we mention at this point that the      
three-block method and the palindrome method are whole-line methods,      
whereas the two-block method is particularly suitable in our context      
since it does work on a half-line only and the additional input of      
transfer matrix trace estimates is most conveniently studied by means of      
a renormalization approach which works for substitution models as well      
as for Sturmian models.\\[0.5cm]      
We will therefore now recall the two-block method in a form suitable to      
our approach. Applications to Sturmian models are presented in the next section. Define       
$$G(n,C) = \{ \omega \in \Omega \; : \; V_\omega(k) = V_\omega(k+n), 1      
\le k \le n, \; |{\rm tr}(M(E,\omega_,1,n))| \le C \;\; \forall E      
\in \Sigma \}.$$       
      
\begin{prop}\label{method}      
Suppose there is $C < \infty $ such that $\limsup_{n \rightarrow \infty}      
\mu(G(n,C)) > 0$. Then, $\mu (\Omega^s) >0$.       
\end{prop}      
{\it Proof.} By assumption, we have $\mu( \limsup (G(n,C))) > 0$. The assertion follows from (\ref{omegas}) if we show      
\begin{equation}\label{stabcrit}      
\limsup (G(n,C)) \subseteq \Omega^s_+.      
\end{equation}      
Consider an arbitrary $\omega_0 \in \limsup (G(n,C))$. Then, for a      
suitable sequence $n_k \rightarrow \infty$, $\omega_0$ belongs to      
$G(n_k,C)$. For every $k\in \NN$, we therefore have    
$M(E,\omega_0,1,2n_k)=M(E,\omega_0,1,n_k)^2$ and      
$|{\rm tr}(M(E,\omega_0,1,n_k))| \le C$ for every $E \in \Sigma$. Together with the identity      
$$M(E,\omega_0,1,n_k)^2 - {\rm tr}(M(E,\omega_0,1,n_k)) M(E,\omega_0,1,n_k)+ I = 0,$$      
valid for all matrices with determinant 1, this implies for every solution $u$ of (\ref{eve}),      
$$\left(\begin{array}{c} u(2n_k + 1) \\ u(2n_k) \end{array}\right) -      
{\rm tr}(M(E,\omega_0,1,n_k)) \left(\begin{array}{c} u(n_k + 1) \\ u(n_k)      
  \end{array}\right) + \left(\begin{array}{c} u(1) \\ u(0)      
  \end{array}\right) = 0.$$      
Using now the bound on the trace, $|{\rm tr}(M(E,\omega_0,1,n_k))| \le C$, we see that $u$ is not square-summable at $+\infty$ (and in fact not even  converging to zero). Thus, we obtain $\omega_0 \in \Omega^s_+$ and hence the assertion  (\ref{stabcrit}).\hfill$\Box$    
\\[5mm]     
We close this section with a brief discussion of the effects that local perturbations of the operators under study may cause. Recall that apart from its connection with zero-measure spectrum, the definition of the stability set $\Omega^s$ was motivated by the question whether local perturbations can introduce eigenvalues inside the spectrum of the unperturbed operator. Under the assumption of Theorem \ref{basic} (i.e., positivity of $\mu (\Omega^s)$), it follows immediately that for almost every $\omega \in \Omega$ and every finitely supported perturbation $W : \ZZ \rightarrow \RR$, the operator $H_\omega + W$ has no eigenvalues in $\Sigma$. Note, however, that this does not imply $\sigma_{{\rm pp}}(H_\omega + W) \cap \Sigma = \emptyset$. In fact, it can very well happen that $H_\omega + W$ has pure point spectrum. To illustrate this phenomenon let us consider energies in $\RR \setminus \Sigma$. By general principles, $\gamma(E) >0$ for every such energy \cite{cl}. In particular, as soon as the spectral measure of a perturbed operator $H_\omega + W$ has some weight outside of $\Sigma$, the spectrum there must be pure point (since it is discrete) and the corresponding eigenfunctions decay exponentially at the rate of the Lyapunov exponent if $\gamma^\pm$ happen to exist at these eigenvalues.     
\\[5mm]    
More drastically, using spectral averaging methods, one can prove the following localization result. For $\omega \in \Omega$ and $\alpha,\beta \in \RR$, define the perturbed operators $H_{\omega,\alpha,\beta}$ by   
$$   
H_{\omega,\alpha,\beta} = H_\omega + \alpha \langle \delta_0, \cdot   
\rangle \delta_0 + \beta \langle \delta_1, \cdot \rangle \delta_1.   
$$   
   
\begin{prop}   
If $|\Sigma| = 0$, then for almost every $\omega$ and almost every pair $(\alpha,\beta)$, $H_{\omega,\alpha,\beta}$ has pure point spectrum with exponentially decaying eigenfunctions.   
\end{prop}   
{\it Proof.} For $\omega \in \Omega$, let $N_\omega = \{ E \in \RR : \omega \not\in \Omega_E \}.$ Then for almost every $\omega \in \Omega$, we have  $\sigma(H_\omega) = \Sigma$ and $|N_\omega| = 0$. Consider for such an $\omega$ and some fixed $\beta \in \RR$ the one-parameter family $H_\alpha = H_{\omega,\alpha,\beta}$. Then, by a standard spectral averaging argument \cite{s1,s2,sw}, for almost every $\alpha$, the spectral measure of the pair $(H_{\omega,\alpha,\beta},\delta_0)$ has zero weight on $\Sigma \cup N_\omega$. In particular, this spectral measure is concentrated on a countable set. Similarly, we have that for every $\alpha$ and almost every $\beta$, the spectral measure of the pair $(H_{\omega,\alpha,\beta},\delta_1)$ has no weight on $\Sigma \cup N_\omega$ and thus is pure point. Hence for almost every pair $(\alpha,\beta)$, we have that the spectral measures of the pairs $(H_{\omega,\alpha,\beta},\delta_0)$,   
$(H_{\omega,\alpha,\beta},\delta_1)$ are pure point measures supported on the complement of $\Sigma \cup N_\omega$. Since $(\delta_0,\delta_1)$ are cyclic for $H_{\omega,\alpha,\beta}$, this operator has pure point spectrum. Moreover, the corresponding eigenfunctions decay exponentially at the rate given by the Lyapunov exponent (recall that on $\Sigma^c$, $\gamma$ is always strictly positive, and note that local perturbations do not affect the existence of solution asymptotics).\hfill$\Box$     
\\[5mm]    
The proposition applies in particular if the family $(H_\omega)_{\omega \in \Omega}$ obeys the assumption of Theorem \ref{basic}. Note that the proposition and its proof are actually valid for all $\omega\in \Omega$ rather than for almost all $\omega$ if the Lyapunov exponent exists uniformly on $\RR \setminus \Sigma$ and $\sigma(H_\omega)=\Sigma$ for all $\omega \in \Omega$. This holds for instance for models generated by primitive substitutions \cite{h} or Sturmian sequences where the rotation number has bounded continued fraction expansion coefficients \cite{dl2}.  
\\[5mm]  
To summarize, by perturbing the potential at two consecutive sites one can shift the entire mass of the spectral measure from $\Sigma$ to $\RR \setminus \Sigma$, thereby forcing the spectral type to become pure point, and the decay of the corresponding eigenfunctions can be investigated by means of the Lyapunov exponent. However, since local perturbations leave the essential spectrum invariant, we infer that $\Sigma$ is precisely the set of accumulation points of the eigenvalues in the gaps of $\Sigma$. In particular, in the case where $(H_\omega)_{\omega \in \Omega}$ obeys the assumption of Theorem \ref{basic}, for almost every $\omega$, we have $\sigma_{{\rm pp}}(H_\omega) \cap \Sigma = \Sigma$, although $\Sigma$ does not contain an eigenvalue of $H_\omega$.

\section{Application to Sturmian potentials}      
In this section we study the stability set corresponding to families of discrete one-dimensional Schr\"odinger operators with Sturmian potentials. In the appendix we shall show that circle maps give rise to two different dynamical systems (on the torus and on the set of pointwise accumulation points), both being strictly ergodic, which, however, are isomorphic as measure-preserving systems. We shall prove the following theorem, where we consider the hull $\Omega_\alpha$ as defined via accumulation points and the measure given on cylinder sets by frequencies of words. Let us remark that in torus representation, this result can be obtained likewise by a direct proof, using \cite{k1}, without making reference to the correspondence given in the appendix.      
      
\begin{theorem}\label{sturm}      
Let $\alpha \in \TT$ be irrational and consider the Sturmian hull $\Omega_\alpha$. Then, for every $\lambda$, the stability set $\Omega_\alpha^s$ has positive measure.      
\end{theorem}      
Before giving the proof of Theorem \ref{sturm}, let us recall some basic properties of Sturmian hulls. Consider the continued fraction expansion of $\alpha$,      
      
\begin{equation}      
  \label{continuedfraction}      
  \alpha = \cfrac{1}{a_1+ \cfrac{1}{a_2+ \cfrac{1}{a_3 + \cdots}}}      
  = [a_1,a_2,a_3,...]      
\end{equation}      
with uniquely determined $a_n \in \NN$. The associated rational      
approximants $\frac{p_n}{q_n}$ are defined by the following recursions,       
$p_0 = 0, \; p_1 = 1, \; p_n = a_n p_{n-1} + p_{n-2}$ and      
$q_0 = 1, \; q_1 = a_1, \; q_n = a_n q_{n-1} + q_{n-2}$. Define the words $s_n$ over the alphabet $\{0,1\}$ by      
      
\begin{equation}      
  \label{recursive}      
s_{-1} = 1, \;\; s_0 = 0, \;\; s_1 = s_0^{a_1 - 1} s_{-1},      
\;\; s_n = s_{n-1}^{a_n} s_{n-2}, \; n \ge 2.      
\end{equation}      
In particular, the word $s_n$ has length $q_n$, $n \ge 0$.      
By definition, for $n \ge 2$, $s_{n-1}$ is a prefix of $s_n$. Therefore,      
we may define a one-sided infinite sequence in the following way,      
      
\begin{equation}      
  \label{standard}      
c_\alpha = \lim_{n \rightarrow \infty} s_n.        
\end{equation}      
The following proposition connects $c_\alpha$, and thus the words $s_n$, to the Sturmian hull $\Omega_\alpha$.      
      
\begin{prop}\label{vsymm}      
$v_{\alpha}$ restricted to $\{1,2,3,\ldots\}$ coincides with $c_\alpha$.      
\end{prop}      
{\it Proof.} This was shown by Bellissard et al. in      
\cite{bist}.\hfill$\Box$\\[0.5cm]      
The following elementary formula is of profound use, as will be seen below.      
      
\begin{prop}\label{wunderformel}      
For each $n \ge 2$, $s_n s_{n+1}= s_{n+1} s_{n-1}^{a_n - 1} s_{n-2} s_{n-1}$.      
\end{prop}      
{\it Proof.} $s_n s_{n+1} = s_n s_n^{a_{n+1}} s_{n-1} = s_n^{a_{n+1}}      
s_n s_{n-1} = s_n^{a_{n+1}} s_{n-1}^{a_n} s_{n-2} s_{n-1} = s_{n+1}      
s_{n-1}^{a_n - 1} s_{n-2} s_{n-1}.$\hfill$\Box$\\[0.5cm]      
Write $s_n = s_n^1 \ldots s_n^{q_n}$ and define the transfer matrices $M_{\lambda,\alpha,E}(n)$ corresponding to the words $s_n$ by      
$$M_{\lambda,\alpha,E}(n) =       
\left(       
  \begin{array}{cc}      
E-\lambda s_n^{q_n} & -1\\1 & 0      
  \end{array}      
\right)      
\times \cdots \times      
\left(       
  \begin{array}{cc}      
E-\lambda s_n^1 & -1\\1 & 0      
  \end{array}      
\right).$$      
      
\begin{prop}\label{tracebound}      
For every $\lambda \not= 0$, there exists $C_\lambda$ such that, for      
every irrational $\alpha$, every $E \in \Sigma$ and      
every $n \in \NN$, we have $|{\rm tr}(M_{\lambda,\alpha,E}(n)| < C_\lambda$.      
\end{prop}      
{\it Remark.} The dependence of $\Sigma$ on $\lambda$ and $\alpha$ is left implicit.\\[0.5cm]      
{\it Proof.} See \cite{bist}.\hfill$\Box$\\[0.5cm]      
We are now in position to give the proof of Theorem \ref{sturm}.\\[0.5cm]      
{\it Proof of Theorem \ref{sturm}.} We shall rely on an application of Proposition \ref{method}. Since Proposition \ref{tracebound} provides us with a uniform upper bound on the traces of transfer matrices corresponding to occurrences of $s_n$ or cyclic permutations thereof, we only need to estimate the measure of certain cylinder sets from below. Namely, by       
      
$$G(q_n,C_\lambda) = \{ \omega \in \Omega_\alpha : V_\omega(k) = V_\omega(k+n), 1 \le k \le q_n, |{\rm tr}(M(E,\omega_,1,q_n))| \le C_\lambda \; \forall E \in \Sigma \} \supseteq \bigcup_{c_n \sim s_n}[c_n c_n]_{[1,2q_n]}.$$      
Thus, we obtain the bound (compare, e.g., \cite{d3})      
      
\begin{equation}\label{lower}      
\mu(G(q_n,C_\lambda)) \ge q_n \cdot d(s_n s_n s_n).      
\end{equation}      
It therefore remains to be shown that       
      
\begin{equation}\label{freqb}      
\limsup_{n \rightarrow \infty} q_n \cdot d(s_n s_n s_n) > 0.      
\end{equation}      
By Proposition \ref{vsymm} we can investigate the frequency of the word $s_n s_n s_n$ by looking at $c_\alpha$. This is convenient since the local structure of $c_\alpha$ is entirely determined by (\ref{recursive}) and (\ref{standard}). Recall that the word $s_n$ has length $q_n$. The bound (\ref{freqb}) follows once we show that, for infinitely many $n$ and every $k \in \NN$, the word $c_\alpha(k) \ldots c_\alpha(k + 7q_n -1)$ contains a copy of the word $s_n s_n s_n$. We shall refer to this property as the window property. Let us first consider the case $\limsup_{n \rightarrow \infty} a_n \ge 2$. Then, for every $n \ge 3$ with $a_{n+1} \ge 2$, the window property follows from (\ref{recursive}) and (\ref{standard}) along with Proposition \ref{wunderformel}. In fact, in this case a window length of $6 q_n$ will already suffice. If $\limsup a_n$ equals $1$, pick $n_0 \ge 3$ such that $a_n = 1$ for every $n \ge n_0$. For these $n$, the window property with window length $7q_n$ can be established by the same reasoning. Thus in every case, we have (\ref{freqb}). The assertion of the theorem now follows from Proposition \ref{method}.\hfill$\Box$

\section{Extensions, remarks, and open problems}      
In this section we discuss possible extensions of Theorem \ref{sturm}, add some remarks concerning the method employed in comparison to related methods, and list some open problems.\\[0.5cm]      
In the preceding section we have shown that the assumption of Theorem \ref{basic} is satisfied for all Sturmian models. In particular, this yields purely singular continuous spectrum of zero Lebesgue measure. In fact, it can be shown that the stability set equals the entire hull.      
      
\begin{theorem}\label{sturmunif}      
Let $\Omega_\alpha$ be Sturmian. Then, $\Omega_\alpha^s = \Omega_\alpha$.       
\end{theorem}      
A proof of this generalization of Theorem \ref{sturm} can be achieved by the approach developed in \cite{dl1} and \cite{dl4}. It can even be shown that $\Omega_{\alpha,+}^s=\Omega_{\alpha,-}^s=\Omega_\alpha$, so there is no $l^2$-decay of solutions whatsoever (cf. \cite{dl4} for details).    
\\[0.5cm]      
We now turn to a discussion of the method we employed above and compare it to related approaches. Spectral properties of one-dimensional Schr\"odinger operators with strictly ergodic potentials taking finitely many values are intrinsically connected to combinatorial properties of finite and infinite words, the latter being a discipline of its own independent interest (see in particular the monographs \cite{loth1,loth2}). This is basically because of equation (\ref{muprop}) which connects the unique ergodic measure to the frequencies of finite words in infinite words. It is therefore natural to look for connections between combinatorial properties of infinite words and spectral theory when studying the spectral properties of the associated operators $H_\omega$. In this regard, basically two fundamental connections have been found. Gordon proved in 1976 that powers occurring in infinite words enable one to deduce bounds on generalized eigenfunctions \cite{g}, allowing for the exclusion of eigenvalues. While his original result essentially requires the occurrence of a fourth power (see also \cite{cfks}), it was later shown by Delyon-Petritis \cite{dp1} and S\"ut\H{o} \cite{s3} that cubes and squares lead to useful estimates, too. We shall refer to these variants of the Gordon criterion as the three-block method and the two-block method, respectively. The two-block method is actually at the heart of Proposition \ref{method}. Among the works that have applied a Gordon-type argument, we want to mention \cite{bist,bg1,bg2,d1,d2,d3,dl4,dl1,dp1,k1,s3}. On the other hand, based on Jitomirskaya-Simon \cite{js}, Hof et al. have shown that palindromes occurring in infinite words also imply certain restrictions on generalized eigenfunctions \cite{hks}. We shall refer to their criterion as the palindrome method.\\[0.5cm]      
What all these methods have in common is a certain kind of limiting procedure which takes finite local symmetries (such as powers or palindromes) to an appropriate property of an infinite word. This is of course non-trivial, as can be seen already by the fact that in general, strong approximation is not a very good type of convergence when studying the point spectrum. The methods therefore have to conserve a certain quantity when taking the limit, and this quantity is always a uniform lower bound on some part of the generalized eigenfunction. This explains why all these methods (with some restrictions in case of the palindrome method) not only exclude eigenvalues but also the existence of (globally) decaying solutions.\\[0.5cm]      
This brings us to our central objective here. We have seen that the    
non-existence of decaying solutions is connected to the zero-measure property. However, this non-existence additionally needs some uniformity in the choice of half-line (at least for fixed energy) having no decaying solution. This immediately rules out the three-block method and the palindrome method since they do not give any information as to where in terms of half-lines the part of the generalized eigenfunction which can uniformly be estimated is located. This leaves us with the two-block method. We have seen that this method indeed allows one  to provide a  criterion  for the    
stability set not to  be empty and that this criterion  is applicable to all Sturmian models.\\[0.5cm]      
We close this section with a list of open problems we consider interesting in this context.      
      
\begin{enumerate}      
\item Prove that, in the non-Sturmian circle map case, the spectrum has zero Lebesgue measure. Essentially, this amounts to finding a suitable analog to the trace map in order to relate the spectra of periodic approximants to the spectrum of the limiting operator, but one could also imagine other ways to establish a vanishing Lyapunov exponent for energies from the spectrum.       
\item Study the quantum dynamics, particularly those of Sturmian models. Partial results in this direction are contained in \cite{d1,dl4,jl1,kkl}.      
\item Identify the Hausdorff dimension of both the the spectrum and the spectral measures. Lower bounds were obtained in \cite{d1,dl4,jl1} for a zero-measure set of $\alpha$'s, all $\lambda$ and all elements in the hull; the only known upper bound holds in the Fibonacci case for large coupling $\lambda$ \cite{r}. Thus, for many parameter values essentially nothing is known rigorously about the dimensional properties of the spectrum and the spectral measures. But already the known bounds may not be optimal. Hence there is room for many new ideas in this direction.      
\item Find an example of a substitution or circle map model with non-empty point spectrum.          
\end{enumerate}

\section{Appendix: Strictly ergodic dynamical systems associated to      
  circle maps, torus vs. accumulation points}      
Consider the family $H_\theta$, $\theta \in [0,1)$, of  discrete one-dimensional Schr\"odinger operators with circle map potentials, that is,      
$$H_\theta u (n)= u(n+1) + u(n-1) + v_\theta(n) u(n),$$      
where $v_\theta(n)= \lambda \chi_{[1-\beta, 1)}(\alpha n + \theta \mod 1)$ with $\lambda \not= 0$, $\alpha \in (0,1)$ irrational and $\beta \in (0,1)$ arbitrary. When studying this family of operators, one tries to fit it  into the framework of random Schr\"odinger operators. Thus, one tries to provide a dynamical system $(\Omega, \mu, T)$ in the sense of the introduction of this paper such that the family $(H_\theta)_{\theta \in [0,1)}$ is part of the family $(H_\omega)$. Quite surprisingly, it turns out that such a system is not uniquely determined, even if one requires the system to be minimal and the family of operators to have the property that the orbit of each operator is dense with respect to strong convergence. In fact, two different dynamical systems and consequently two different families of operators have been investigated in the past. It is the aim of this appendix to study the mutual relation between these systems. After some preliminary propositions we give the connection between the underlying sets in Lemma \ref{set}. In Theorem \ref{homeo} we provide a homeomorphism between (rather large) invariant subsystems of the respective systems. Finally, Corollary \ref{iso} shows that the two systems are isomorphic from a measure theoretical point of view.\\[0.5cm]      
Let an irrational $\alpha\in (0,1)$ and an arbitrary $\beta\in (0,1)$ be      
given. We start by recalling the construction of the two dynamical      
systems in question, where, for simplicity, we restrict our attention to the case $\lambda =1$. The general case can of course be treated similarly.        
      
\begin{itemize}      
\item Identify $I:=[0,1)$ with $\RR/\ZZ$. Then $I$ is compact and  metrizable. Consider the shift      
$$S: I\longrightarrow I, \;\: x\mapsto x+\alpha \mod 1.$$      
Then $(I,S)$ is a strictly ergodic dynamical system \cite{w}. Denote the      
unique normed invariant measure on $I$ by $\mu_I$. This system gives rise to a family of Schr\"odinger operators $(H_\theta)$ given by      
$$(H_\theta \phi)(n) = \phi (n+1) + \phi (n-1) + v_{\theta}(n) \phi      
(n),$$      
with $v_\theta(n)=\chi_{[1-\beta,1)}(\alpha n + \theta \mod 1)$.      
This family of operators has been studied in \cite{bist,dp1,k1}.      
\item Define $\Omega\subset \{0,1\}^{\ZZ}$ by       
$$\Omega=\{ \lim T^{n_i} v_0\,:\, n_i\to \infty\},$$      
where $T:\{0,1\}^{\ZZ} \rightarrow \{0,1\}^{\ZZ}$, $Tu(n)=      
u(n+1)$, denotes the shift. Then $\Omega$ is a compact and metrizable      
space. Again $(\Omega,T)$ is a strictly ergodic dynamical system      
\cite{h}, whose unique normalized invariant measure will be denoted by      
$\mu_\Omega$.  The corresponding family of Schr\"odinger operators is      
given by       
$$(H_\omega \phi)(n) = \phi (n+1) + \phi (n-1) + V_\omega (n) \phi (n),$$      
where $V_\omega(n)=\omega(n)$. The investigations in \cite{h,hks} are based upon this system.      
\end{itemize}      
{\it Remark.} In fact, our method can easily be adapted to treat      
functions $v_0$ of the form      
$$v_0(n)=\chi_{I-A}(\alpha n \mod 1),$$      
where $A$ is a finite union of half-open intervals of the form      
$[a,b)$.\\[0.5cm]      
The analysis of the connection between these two systems will be based      
on the study of the continuity properties of the maps      
$$ j:I\longrightarrow  \{0,1\}^{\ZZ}, \; j(\theta) = v_\theta$$      
and      
$$ j_n:I\longrightarrow \RR, \; j_n(\theta)= v_\theta(n),$$      
where $n$ is an arbitrary element of $\ZZ$. We will need some notation:      
Let $(\theta_n)$ be a sequence in $I$ converging to $\theta$ in $I$.  We      
say that the sequence $(\theta_n)$ converges to $\theta$ from the right      
written as       
$$\theta_n \longrightarrow \theta +, \, n\to \infty,$$       
if there exist $\delta>0$ with $[\theta, \theta + \delta) \subset I$      
and $n_0\in \NN$ such that $\theta_n\in [\theta,\theta + \delta)$ for all      
$n\in \NN$ with $n\geq n_0$. Similarly, we say that the sequence  $(\theta_n)$ converges to $\theta$ from the left written as       
$$\theta_n \longrightarrow \theta - , \, n\to \infty,$$       
if there exist $\delta>0$ with $[\theta, \theta + \delta)\subset I$      
and $n_0\in \NN$ such that $\theta_n\notin  [\theta,\theta + \delta)$ for all $n\in \NN$ with $n\geq n_0$. Define for $\theta \in I$ the orbit $O(\theta)$ by       
$$O(\theta)=\{S^n(\theta)\,:\,n\in \ZZ\}=\{\theta +\alpha n\,:\,n\in      
\ZZ\}$$      
and define for $\omega\in \Omega$ the orbit $O(\omega)$ by      
$$O(\omega)= \{T^n(\omega)\,:\,n\in \ZZ\}.$$      
The following proposition is elementary but crucial.       
      
\begin{prop} \label{continuityproperties}      
{\rm (i)} Let $n\in \ZZ$ be given.  The map $j_n$ is continuous at all $\theta\in I$ obeying $\theta + \alpha n \notin \{0, 1-\beta\}$.\\      
{\rm (ii)} If $\theta + \alpha n \in \{0, 1-\beta\}$, then       
$$\lim_{\theta' \to \theta +} j_n(\theta')= j_n(\theta)=\left\{       
  \begin{array}{r@{\quad:\quad}l}      
0 &  \theta + \alpha n = 0 \\      
1 &  \theta + \alpha n = 1-\beta      
  \end{array}      
\right.       
\:\;\mbox{and}\;\:       
\lim_{\theta'\to \theta-} j_n(\theta)= \left\{       
  \begin{array}{r@{\quad:\quad}l}      
1 &  \theta + \alpha n = 0 \\      
0 &  \theta + \alpha n = 1-\beta      
  \end{array}      
\right. . $$      
{\rm (iii)} Let $\theta \in I$ be given. Then there are at most two different $n\in \ZZ$ such that $j_n$ is not continuous at $\theta$. In particular there exists $n_1(\theta)\in \ZZ$ such that $j_n$ is continuous at $\theta$ for all $n \geq n_1(\theta)$.\\      
{\rm (iv)} The map $j$ is continuous at all points $\theta\in I$ with $\theta\notin O(0)\cup O(1-\beta).$      
\end{prop}      
{\it Proof.} A direct calculation yields (i) and (ii). Moreover,      
(iii) and (iv) follow easily from (i).\hfill $\Box$      
      
\begin{prop} \label{o1} The limits       
$$\omega_0=\lim_{\theta\to 0-}      
j(\theta)\;\:\mbox{and}\;\:\omega_{(1-\beta)}=\lim_{\theta\to (1-\beta)-} j(\theta) $$       
exist in $\{0,1\}^{\ZZ}$.  The equations $\omega_0(0)=1$ and      
$\omega_{(1-\beta)}(0)=0$ hold. This implies, in particular, $\omega_0\neq v_0$ and $\omega_{(1-\beta)}\neq v_{(1-\beta)}$. Moreover, there exists $n_2\in \ZZ$ with $\omega_0(n)= v_0(n)$ and $\omega_{(1-\beta)}(n)=      
v_{(1-\beta)}(n)$ for all $n \geq n_2$.      
\end{prop}      
{\it Proof.} This is immediate from Proposition \ref{continuityproperties}.\hfill $\Box$\\[5mm]      
{\it Remark.} In fact, it is not hard to see that      
$\omega_\gamma = \lim_{\theta\to \gamma-} j(\theta)$ exists for all      
$\gamma \in I$. However, for $\gamma \notin O(0)\cup O(1-\beta)$, we have  $\omega_\gamma= v_\gamma$, whereas for $\gamma \in O(0)\cup O(1-\beta)$,      
it is easy to see that $\omega_\gamma$ belongs to $O(\omega_{0})\cup O(\omega_{(1-\beta)})$.      
      
\begin{prop}\label{injectivity}      
Let $\theta$ and $\theta'$ in $I$ be given. If there exists $n_0\in \NN$ with $v_\theta(n)=v_{\theta'}(n)$ for all $n\geq n_0$, then $\theta=\theta'$.      
\end{prop}      
{\it Proof.} This follows easily from the fact that $\{\alpha n + \theta \mod 1\,:\, n\in \NN, \,n \geq n_0\}$ is dense in $I$.\hfill$\Box$\\[5mm]      
We can now describe the relation between $I$ and $\Omega$.      
      
\begin{lemma}\label{set}      
The set $\Omega$ is the disjoint union of $j(I)$ and $O(\omega_{0})\cup O( \omega_{(1-\beta)})$.      
\end{lemma}      
{\it Proof.} We first show  $j(I)\cap (O(\omega_{0})\cup O( \omega_{(1-\beta)}))=\emptyset$: As $j(I)$ is shift invariant, it is enough to show that neither $\omega_0$ nor $\omega_{(1-\beta)}$  belongs to $j(I)$. Assume the contrary. This means that there exists $\theta\in I$ with $\omega_0 =v_\theta$ or $\omega_{(1-\beta)} = v_\theta $. Consider the case $\omega_0 =v_\theta$; the other case can be treated similarly. As $\omega_0(n) = v_0(n)$ for $n\geq n_2$ by Proposition \ref{o1}, we can conclude from Proposition \ref{injectivity} that $\theta=0$. This implies $\omega_0=v_0$. But this is a contradiction by Proposition \ref{o1}.\\[1mm]      
Next, we show $\Omega = j(I)\cup( O(\omega_{0})\cup O( \omega_{(1-\beta)})  )$: The inclusion $j(I)\cup(  O(\omega_{0})\cup O( \omega_{(1-\beta)}) )  \subset \Omega$ is clear. To prove the opposite inclusion, fix $\omega\in \Omega$. By the definition of $\Omega$ we have $\omega=\lim T^{n_i} v_0$. This means that there exists a sequence $(\theta_n)\subset I$ with $\omega=\lim j(\theta_n)$. As $I$ is compact, we can assume w.l.o.g.~that $(\theta_n)$ converges to $\theta\in I$. We will consider three cases:\\[1mm]      
{\it Case 1:} $\theta\notin O(0)\cup O(1-\beta)$. By Proposition \ref{continuityproperties} the map $j$ is continuous at $\theta$. This implies $\omega=\lim j(\theta_n)= j(\theta)\in j(I).$\\[1mm]      
{\it Case 2:} $\theta\in O(0)$. Assume w.l.o.g. $\theta=0$.      
As $j(\beta)(0)=v_\beta(0)=0$ for $\beta< 1-\beta$ and $j(\beta)(0)=v_\beta(0)=1$ for $\beta\geq 1-\beta$ and as the sequence  $(j(\theta_n))$ converges, there are two subcases.\\      
{\it Subcase 1:} $\theta_n< 1-\beta$ for large $n\in \NN$. By Proposition \ref{continuityproperties} the sequence $(j(\theta_n))$ converges to      
$v_0\in j(I)$. This implies $\omega=\lim j(\theta_n)= v_0$.\\      
{\it Subcase 2:} $\theta_n\geq1-\beta$ for large $n\in \NN$. By   Proposition \ref{o1} the sequence  $(j(\theta_n))$ converges to $\omega_0$. This implies $\omega=\lim j(\theta_n)= \omega_0$.\\[1mm]      
{\it Case 3:} $\theta\in O(1-\beta)$. As in Case 2 it is possible to show that either $\omega =T^n v_{(1-\beta)}$ or $\omega =T^n \omega_{(1-\beta)}$ with a suitable $n \in \ZZ$.\\[1mm]      
Thus in all cases, the element $\omega$ belongs to $j(I) \cup      
O(\omega_0)\cup O(\omega_{(1-\beta)})$, concluding the proof.\hfill $\Box$      
      
\begin{theorem}\label{homeo}      
The restriction $\tilde{j}$ of $j$ to $\tilde{I}=  I- (O(0)\cup O(1-\beta)) $ is a homeomorphism between $\tilde{I}$  and $\tilde{\Omega}= \Omega-(O(v_0)\cup O(\omega_0)\cup O(v_{(1-\beta)})\cup O(\omega_{(1-\beta)})  )$.      
\end{theorem}      
{\it Proof.} Using Proposition \ref{injectivity}, it is not hard to show      
that $j$ maps $\tilde{I}$ into $\Omega-(O(v_0) \cup O(v_{(1-\beta)}))$. By Lemma \ref{set} the image of $I$ under $j$ is contained in $\Omega- (O(\omega_0) \cup O(\omega_{(1-\beta)}))$. Thus $j$ is indeed a map between      
$\tilde{I}$ and $\tilde{\Omega}$. It is one to one by Proposition \ref{injectivity} and onto by Lemma \ref{set}. It is continuous by Proposition \ref{continuityproperties}.\\[1mm]      
It remains to be shown that its inverse is  continuous. Suppose $(\omega_n)$ is a sequence in $\tilde{\Omega}$ converging to $\omega\in  \tilde{\Omega}$. By the already proven bijectivity  of $j$, there exists a sequence $(\theta_n)\subset I$ with $\omega_n=j(\theta_n) $. As already shown, the map $j$ is continuous at every $\theta \in \tilde{I}$. Thus, we have to show that $(\theta_n)$ converges to $\theta \in \tilde{I}$. By $(\theta_n)\subset I$ and the compactness of $I$ there exists a subsequence of $(\theta_n)$ converging in $I$. It is not hard to see, using Proposition \ref{continuityproperties} and Proposition \ref{injectivity}, that all converging subsequences must have the same limit in $I$. This implies that the sequence $(\theta_n)$ does indeed converge to $\theta \in I$. We finish the proof by showing that $\theta$ does not belong to $O(0)\cup O(1-\beta)$. Assume the contrary.  As in the proof of Case 2 and Case 3 of Lemma \ref{set}, it is then possible to show that $\omega=\lim j(\theta_n)$  belongs to $(O(v_0)\cup O(\omega_0)\cup O(v_{(1-\beta)})\cup O(\omega_{(1-\beta)}))$. This is a contradiction as $\omega$ belongs to $  \tilde{\Omega}$.\hfill $\Box$\\[5mm]      
{\it Remark.} Using Theorem \ref{homeo} and the fact that $(I,S)$ is minimal, it is possible to show that $(\Omega,T)$ is minimal. This gives an      
alternative to the proof in \cite{h}.\\[0.5cm]      
Recall that two measure-preserving systems $(M_i,\mu_i,S_i)$, $i=1,2$,       
where $M_i$ is a measure space, $\mu_i$ is a measure on $M_i$, and $S_i$ is      
a measure-preserving map on $M_i$, are called isomorphic if there      
exist invariant subsets $\tilde{M_i}$ of $M_i$ having full measure and a      
bimeasurable, measure-preserving map $\phi: \tilde{M_1}\longrightarrow \tilde{M_2}$ with      
$$\phi\circ S_1 = S_2 \circ \phi$$      
on $\tilde{M_1}$ (cf. Definition 2.4 in \cite{w}).       
      
\begin{coro}\label{iso}      
The measure-preserving systems $(I, \mu_I, S)$ and $(\Omega, \mu_\Omega, T)$ are isomorphic.      
\end{coro}      
{\it Proof.} Obviously, the set $O(0)\cup O(1-\beta)$ (resp., $O(v_0)\cup  O(\omega_0)\cup O(v_{(1-\beta)})\cup O(\omega_{(1-\beta)})$) is invariant in $I$ (resp., $\Omega$). Moreover, it has measure zero, as an orbit consisting      
of infinitely many points in a space of finite invariant measure has measure zero.      
Thus, Theorem \ref{homeo} immediately implies the existence of invariant sets $\tilde{I}$ and $\tilde{\Omega}$ of full measure and the existence of a bijective, bimeasurable map $\tilde{j}: \tilde{I}\longrightarrow \tilde{\Omega}$. A short calculation gives $\tilde{j} \circ S = T \circ \tilde{j}$. This equation, the uniqueness of the invariant measure on $I$ (resp., $\Omega$), and the bimeasurability of $\tilde{j}$ imply that $\tilde{j}$ is measure preserving. \hfill $\Box$\\[5mm]      
{\it Remark.} By a slight change of the argument in the corollary, it      
is possible to conclude the unique ergodicity of the system $(\Omega,T)$      
from the unique ergodicity of the system $(I,S)$.       
\\[5mm]      
{\it Acknowledgments.} The authors would like to thank R.~del Rio for a      
discussion initiating this study. D.~D.~was supported by the German Academic Exchange Service through Hochschulsonderprogramm III (Postdoktoranden) and D.~L.~received financial support from Studienstiftung des Deutschen Volkes, both of which are gratefully acknowledged.

\end{document}